\documentclass[12pt]{article} 
\usepackage{fullpage,amsmath,amscd,amsthm,amssymb,amsxtra,latexsym,epsfig, 
epic,eepic,mathrsfs} 
\begin{document} 
\title{\bf Tjurina and Milnor numbers of matrix singularities} 
\author{V.~Goryunov \quad and \quad D.~Mond} 
\date{} \maketitle
                        
\def\RR{{\bf R}}
\def\AA{{\cal A}}
\def\BB{{\cal B}}
\def\KK{{\bf K}}
\def\LL{{\bf L}}
\def\NN{{\bf N}}
\def\ZZ{{\bf Z}}
\def\CC{{\bf C}}
\def\JJ{{\bf J}}
\def\Proof{{\em Proof}\quad}
\def\SS{{\mathscr S}}
\def\K{{\cal K}}
\def\M{{\cal M}}
\def\P{{\cal P}}
\def\OO{{\cal O}}
\def\S{{\Sigma}}
\newcommand\s{\mathscr}
\def\mod{\,{\rm mod}\,}
\def\a{{\alpha}}
\def\b{{\beta}}
\def\f{{\varphi}}
\def\F{{\Phi}}
\def\ni{{\noindent}}
\def\l{{\lambda}}
\def\L{{\Lambda}}
\def\O{{\Omega}}
\def\o{{\omega}}
\def\G{{\Gamma}}
\def\GG{{\cal G}}
\def\g{{\gamma}}
\def\k{{\kappa}}
\def\ben{\begin{enumerate}}
\def\en{\end{enumerate}}
\def\bin{\begin{itemize}}
\def\ein{\end{itemize}}
\def\i{\item}
\def\eop{{\hfill{$\Box$}}\\}
\def\Ext{{\mbox{Ext}}}
\def\Tor{{\mbox{Tor}}}
\def\Der{{\mbox{Der}}}
\def\codim{{\mbox{codim}}}
\def\gl{{\mbox{gl}}}
\def\Gl{{\mbox{Gl}}}
\def\Pf{{\mbox{Pf}}}
\def\sPf{{\mbox{\footnotesize Pf}}}
\def\pf{{\mbox{\footnotesize Pf}}}
\def\sl{{\mbox{sl}}}
\def\Sl{{\mbox{Sl}}}
\def\Mat{{\mbox{Mat}}}
\def\sMat{{\mbox{\footnotesize Mat}}}
\def\mat{{\mbox{mat}}}
\def\m{m_0}
\def\sk{{\mbox{sk}}}
\def\Sk{{\mbox{Sk}}}
\def\sSk{{\mbox{\footnotesize Sk}}}
\def\sym{{\mbox{sym}}}
\def\Sym{{\mbox{Sym}}}
\def\sSym{{\mbox{\footnotesize Sym}}}
\def\JP{{\mbox{\bf JP}}}
\def\J{{\mbox{\bf J}}}
\def\k{{\mbox{\bf K}}}
\def\L{{\mbox{\bf L}}}
\def\F{{\mbox{\bf F}}}
\def\f{{\s F}}
\def\GN{{\mbox{\bf GN}}}
\def\K{{\cal K}}
\def\D{{\Delta}}
\def\dd{{\partial}}
\def\p{{\partial}}
\def\DD{{\cal D}}
\def\d{{\delta}}
\def\t{{\theta}}
\def\T{{\Theta}}
\def\x{{\chi}}
\def\e{{\varepsilon}}
\def\ti{\widetilde}
\def\ld{\ldots}
\def\cd{\cdots}
\def\ha{\widehat}
\def\R{{\cal R}}
\def\C{{\cal C}}
\def\la{{\langle}}
\def\ra{{\rangle}}
\def\dim{{\rm dim}}
\def\be{\begin{equation}}
\def\ee{\end{equation}}
\def\iff{{\Leftrightarrow}}
\def\imp{{\Rightarrow}}
\def\to{{\ \rightarrow\ }}
\def\too{{\ \longrightarrow\ }}
\def\st{\stackrel}
\def\bu{{\bullet}}
\def\bw{\bigwedge}
\def\da{\downarrow}
\newtheorem{theorem}{Theorem}[section]
\newtheorem{definition}[theorem]{Definition}
\newtheorem{remark}[theorem]{Remark} 
\newtheorem{lemma}[theorem]{Lemma}
\newtheorem{example}[theorem]{Example}
\newtheorem{proposition}[theorem]{Proposition}
\newtheorem{corollary}[theorem]{Corollary}
\newtheorem{conjecture}[theorem]{Conjecture}
\newtheorem{remarks}[theorem]{Remarks}

In \cite{B}, Bill Bruce classified simple 
singularities of symmetric matrix families with respect to a 
natural equivalence relation (see Section \ref{mdets} below). 
The very first look at his tables reveals
a rather unexpected relation peculiar to two-parameter
families: the dimension of the base of a matrix miniversal
deformation coincides with the Milnor number of the determinant
of the family. This observation was the main motivation
for the paper \cite{BGZ} where equality of the Tjurina and
Milnor numbers for hypersurface sections of an isolated hypersurface
singularity was proven. That provided a partial explanation for the
symmetric matrix problem, in the case of $2 \times 2$ matrices.

In the present paper, we prove that $\tau=\mu$ 
for two-parameter families of symmetric matrices of any order. We also
prove similar statements for two other closely related matrix classification 
problems: for arbitrary square matrices depending on 3 parameters
(as conjectured in \cite{GZ}), and
for skew-symmetric $2k \times 2k$ matrices in 5 variables. The key to
the proof is to switch from the Koszul complex used in \cite{BGZ}
to an appropriate free resolution of a determinantal or Pfaffian variety. 
Very suprisingly, these resolutions
(introduced in \cite{J},\cite{gn} and \cite{JP}) 
have all been known for many years and use exactly the Lie algebras
involved in the matrix classifications introduced just very recently
in \cite{B,BT,BGZ,GZ,Ha}. 

The coincidence of $\tau$ and $\mu$ in the three classifications
is a particular case of our main result, Theorem \ref{betas}, on the relation
between these two numbers for sections, with isolated singularity, of a 
possibly non-isolated hypersurface singularity $(V,0)\subset(\CC^n,0)$,
obtained by pulling back $V$ by a map $F:(\CC^m,0)\to (\CC^n,0)$:
\begin{equation}\label{eq1}
\tau = \mu - \b_0 + \b_1
\end{equation}
where the $\b_i$ are the Betti numbers of the pull-back of a free
resolution of the jacobian algebra of the hypersurface singularity
(and thus are the ranks of certain $\Tor$ modules).
The main result of \cite{BGZ}, 
on sections of an isolated hypersurface singularity, is also a 
special case. 

The formula (\ref{eq1}) has a topological interpretation when the
singular subspace of $V$ is Cohen-Macaulay. In the case
of families of matrices (symmetric, general, skew) in 2 or 3, 3 or 4, 
and 5 or 6 variables respectively, the right-hand side is in fact 
the rank of the vanishing
homology of $V(\det(S_t))$ (or the Pfaffian $V(\Pf(S_t))$ in the 
skew symmetric case) for a generic perturbation $S_t$ 
of the matrix $S$. This is
the generalised Milnor number, in the sense of the rank of the vanishing
homology of the nearby stable object, appropriate to the problem considered. 
In 2, 3 and 5 variables, $V(\det(S))$,  or $V(\Pf(S))$ in the skew 
case, is smoothed in a deformation of $S$; in dimensions $3, 4$ and $6$ 
it is not. In all these cases, however, one has 
$$\mbox{rank of the vanishing homology}=\mbox{Tjurina number}.$$

Closely related results on sections of free divisors were obtained in 
a number of papers of Jim Damon, e.g. \cite{D589},
\cite{leg}, and in the paper \cite{DM} of Jim Damon and the second author. 
Free divisors in $\CC^n$
are characterised by having a singular subspace which is Cohen-Macaulay of 
dimension $n-2$. Our topological results here concern sections 
with isolated singularity
of hypersurfaces whose singular subspace is Cohen-Macaulay of arbitrary 
dimension.
In \cite{DM} the freeness of the
divisors was important in allowing an easy proof of conservation of
multiplicity (analogous to conservation of the Milnor number in a deformation
of a function singularity), which was the basis for the comparison between
the vanishing homology and the module of first-order deformations. Here
we use the conservation of Milnor number itself. 

\medskip 
The structure of the paper is as follows.

The proof of our main theorem here
is a modification of the arguments used (not very explicitly) 
in \cite{BGZ} to prove
(\ref{eq1}) for sections of
isolated hypersurface singularities. We generalise this in Section \ref{def}
by identifying
the $T^1$ of a section of a singular hypersurface as something rather like
(the homology of) a derived
functor associated to the comparison between the jacobian algebra
of the section and the pull-back of the jacobian algebra of the hypersurface
singularity. This allows us to embed the $T^1$ in a canonical long exact 
sequence, from which (\ref{eq1}) follows in case the section has isolated
singularity.


Section \ref{mdets} deals in detail with various notions of 
equivalence of matrix families,
and Section \ref{alg} describes the commutative algebra which allows
us to relate the various notions of equivalence to one another, and
to prove conservation of multiplicity. Section \ref{topdef}
builds on the material of the previous sections to relate $\tau$ to the 
rank of the vanishing homology, and in Section \ref{CM} we comment on
Cohen-Macaulay properties of the relative $T^1$'s and on
some related results which suggest that the discriminants
in the matrix versal base spaces are free divisors.

We are grateful to Dmitry Rumynin for suggesting we use the cone construction.
\subsection{Notation}
At the urging of Kyoji Saito, in this paper we have harmonised our notation 
with standard notation in algebraic geometry. Given a divisor $V$ in a 
complex space 
$X$, it has been usual in singularity-theory papers 
to denote by $\Der(\log V)$ (or $\mbox{Derlog}(V)$)
the $\OO_X$-module of vector fields on $X$ which
are tangent to $V$ at its smooth points, and by $\O^1(\log V)$ the 
$\OO_X$-module of $1$-forms with logarithmic poles along $V$. Since these two
modules are mutually dual, the conventions of 
algebro-geometric notation would insist that
one of them have a $-\log V$ in parentheses rather than $\log V$. As 
a logarithmic pole is after all a pole, and $k$-forms with a (first-order)
pole along $V$ are denoted $\O^k(V)$, it is clear that it has to be
$\Der(\log V)$ that accepts the minus sign and becomes henceforth 
$\Der(-\log V)$. Similarly, we have replaced the notation $\Der(\log f)$, for
the module of vector fields tangent to all the level sets of a function $f$,
by $\Der(-\log f)$.

\section{Equivalence and deformations of sections of hypersurfaces and 
functions}\label{def}
Consider a pair consisting of a function $f:(\CC^n,0)\to (\CC,0)$ and 
a map $F:(\CC^m,0)\to(\CC^n,0)$. Let $V=f^{-1}(0)$. We seek to describe 
the deformations of $F$ in relation to $f$, with a view to understanding 
the deformations of
$F^{-1}(V)$. Many of our formulae will involve both $\OO_{\CC^m}$
and $\OO_{\CC^n}$; where possible, we will use $\OO$ to
abbreviate $\OO_{\CC^m}$, but {\it never} $\OO_{\CC^n}$.
\begin{definition}
{\em 
Two map-germs $F,F': (\CC^m,0) \to (\CC^n,0)$ are called} $\K_f$-equivalent
{\em
if there exist a diffeomorphisms $\Phi$ of $(\CC^m\times\CC^n,0)$
and $\varphi$ of $(\CC^m,0)$, 
such that
\ben
\i
$\pi_1\circ\Phi=\varphi\circ\pi_1$, (i.e. $\Phi$ lifts $\varphi$),
\i  $f\circ\pi_2\circ\Phi=f\circ\pi_2$ (i.e. $\Phi$ preserves $f$)
\i  $\Phi$ induces a diffeomorphism
$\mbox{graph}(F)\to\mbox{graph}(F')$.
\en
}
\end{definition}
This equivalence was introduced in \cite{DM}. It is closely related to
$\K_V$-equivalence, introduced by Damon in \cite{D},
in which (2) is replaced by\\

\hspace{-.05in}2'. $\Phi$  sends $\CC^m\times\{f=0\}$ to itself.\\ 

If $F$ and $F'$ are $\K_f$-equivalent then $f\circ F$ 
and $f\circ F'$ are right equivalent, and if $F$ and $F'$ are $\K_V$
-equivalent then $f\circ F$ and $f\circ F'$ are contact-equivalent.

The extended tangent space to the $\K_f$-orbit of $F$ is
$$
T_{\K_f}F = tF(\t_{\CC^m}) + F^*(\Der(-\log f))\,,
$$
where 
\bin
\i
$\t_{\CC^m}$ is the space of germs of holomorphic vector fields
on $(\CC^m,0)$, 
\i
$tF:\theta_{\CC^m}\to 
\theta(F):=F^*(\theta_{\CC^n})$
is the sheafification of the derivative $dF$, 
\i
and $\Der(-\log f) \subset \t_{\CC^n}$ is the 
$\OO_{\CC^n}$-module of vector fields annihilating $f$, and $F^*(\Der(-\log f))$
is the $\OO_{\CC^m}$-submodule of $\theta(F)$ generated by the composites
with $F$ of the vector fields in $\Der(-\log f)$. 
\ein

The extended tangent space to the $\K_V$ orbit of $F$ is
$$
T_{\K_f}F = tF(\t_m) + F^*( \Der(-\log V))\,,
$$
where $\Der(-\log V)$ is the $\OO_{\CC^n}$-module of vector fields on $\CC^n$
which are tangent to 
$V$ at its smooth points.
Denote $\theta(F)/T\K_fF$ and $\theta(F)/T\K_VF$ respectively
by $T^1_{\K_f}F$ and $T^1_{\K_V}F.$

Notice that 
$$\frac{\theta(F)}{T\K_VF}\otimes_{\OO_{\CC^m,x}}\CC =\frac{T_{F(x)}\CC^n}
{d_xF(T_x\CC^m)+T^{\log}_{F(x)}V}$$ 
where for any point 
$y\in \CC^n,\ T^{\log}_yV=\{\zeta(y):\zeta\in \Der(-\log V)_y\}$ 
is the {\it logarithmic tangent space} to $V$ at $y$. Thus 
$T^1_{\K_V}F$ measures failure of ``logarithmic transversality''
(or algebraic transversality, in Damon's terminology)
of $F$ to $V$. The geometric interpretation of $T^1_{\K_f}F\otimes\CC$ is
less clear --- see Remark \ref{normalforms}(iii) below --- although 
$T^1_{\K_f}F$ in some sense measures
failure of logarithmic transversality of $F$ to the level sets of $f$.

Both $\K_V$ and $\K_f$ are ``geometric subgroups'' of the group of all
diffeomorphism-germs, and so by Damon's general theory (\cite{geosubgps})
the usual theorems of singularity theory 
apply: finite determinacy, infinitesimal criterion for versality, etc.
In particular $T^1_{\K_f}F$ and $T^1_{\K_V}F$ are
the tangent spaces at $0$ to the (smooth) miniversal base-spaces of 
$F$ for the two equivalences. 

Now we describe another approach to these two deformation theories, which
identifies the two $T^1$'s as something resembling a derived functor. 
By this means we are
able to locate them in long exact
sequences which provide solutions\footnote{Of course this is 
not how the solution was first found! Nevertheless this seems 
to be the most canonical presentation.} 
to the problems that motivated this paper.

Consider the two surjective comparison maps 
$$\frac{\OO_{\CC^m}}{J_{f\circ F}}\to \frac{\OO_{\CC^m}}{F^*(J_f)}$$
and
$$\frac{\OO_{\CC^m}}{(f\circ F)+J_{f\circ F}}\to 
\frac{\OO_{\CC^m}}{F^*((f)+J_f)}.$$
By considering free resolutions of the modules involved, we are going to 
incorporate these maps into long exact sequences also involving 
$T^1_{\K_f}F$ and
$T^1_{\K_V}F$. 

Let $\L_\bu$ and $\ti \L_\bu$ be $\OO_{\CC^n}$-free resolutions of 
$\OO_{\CC^n}/J_f$ and $\OO_{\CC^n}/(f)+J_f$, let $\K_{\bu}(f)$ 
be the Koszul complex on the first-order partials
of $f$, and let $\ti \k_\bu(f)$ be $\k_\bu(f)$ augmented by another
generator in degree 1, mapping onto $f$ in degree $0$ (so that 
$H_0(\ti \k_\bu(f)) =\OO_{\CC^n}/(f)+J_f$). 
By lifting the identity maps on $\OO_{\CC^n}/J_f$ and $\OO_{\CC^n}/(f)+J_f$,
we obtain morphisms of complexes 
$$\k_\bu(f)\to \
\L_\bu$$
and 
$$\ti \k_\bu(f)\to \ti \L_\bu.$$
These complexes and morphisms can be pulled back by $F$ --- 
in other words, tensored over $\OO_{\CC^n}$ with $\OO_{\CC^m}$.
There is a natural morphism 
$$\wedge^\bu tF:\k_\bu(f\circ F)\to 
F^*(\k_\bu(f));$$
and by taking the direct sum of $tF$ with the identity map
$\OO\to\OO$ in degree 
$1$ (recall that we frequently abbreviate 
$\OO_{\CC^m}$ simply to $\OO$)
we also obtain a morphism
$$\wedge^\bu \ti {tF}:\ti \k_\bu(f\circ F)\to 
F^*(\ti \k_\bu(f)).$$ 
By composing these with the pulled-back morphisms mentioned above, we obtain
morphisms of complexes
$$\phi_f:\k_\bu(f\circ F)\to F^*(\L_\bu)$$
and
$$\phi_V:\ti \k_\bu(f\circ F)\to F^*(\ti \L_\bu).$$
Let $C_\bu(\phi_f)$ and $C_\bu(\phi_V)$ be the cones (cf \cite{gelman},
pp. 153-158, but the definition is recalled below)
on these morphisms of complexes. Then our main technical result is
\begin{theorem}\label{main} 
$$T^1_{\K_f}F\simeq H_1(C_\bu(\phi_f))$$
and
$$T^1_{\K_V}F\simeq H_1(C_\bu(\phi_V)).$$
\end{theorem}
\ni {\em Proof}\quad Both statements are straightforward consequences
of the definitions. If 
$\psi:A_\bu\to B_\bu$ is a morphism of
complexes then the cone $C_\bu(\psi)$ has $C_n(\psi)=A_{n-1}\oplus B_n$
and differential taking 
$(a_{n-1},b_n)$ to $(-d(a_{n-1}),d(b_n)-\psi(a_{n-1}))$. 

The following
diagram shows the morphism $\phi_f$.

\be\label{bigdiag}
\begin{array}{clcccccccccccccc}
\k_\bu(f\circ F) & : & \qquad & \dots & \st{d(f\circ F)}\too & \bw^2 \OO^m &
\st{d(f\circ F)}\too & \OO^m & \st{d(f\circ F)}\too & \OO & \too & 0 \\
\\
\phi_f\da  &&&  \dots & & 
\da & &  
\hspace{-15pt}{}_{dF} \da && 
\hspace{-9pt}{}_{id} \da \\
\\
F^*(\L_\bu) &: & & \dots & \st{}\too & \OO^r &
\st{F^*(\alpha_1)}\too & F^*(\theta_{\CC^n}) & \st{F^*(df)}\too & \OO & \too & 0 \\ 
\end{array}
\ee
The cone is the total complex of this (rather small) double 
complex; its modules
are direct sums along the south-west to north-east parallels, and its 
differential runs south-east.
Note that the image of $F^*(\alpha_1)$ is $F^*(\Der(-\log f))$.

We have 
$$Z_1(C_\bu(\phi_f))=\{(a,\xi)\in \OO\oplus F^*(\theta_{\CC^n}):
F^*(tf)(\xi)=a\}=\{(F^*(tf)(\xi), \xi):\xi\in F^*(\theta_{\CC^n})\},$$
and $H_1(C_\bu(\phi_f))$ is the quotient of this by 
$$\{(-t(f\circ F)(\eta), \zeta-tF(\eta)):\eta\in \theta_{\CC^m}, \zeta\in 
F^*(\Der(-\log f))\}.$$
Under projection to $F^*(\theta_{\CC^n})$ forgetting the first component, 
this quotient maps isomorphically
to 
$$\frac{F^*(\theta_{\CC^n})}{\{\zeta-tF(\eta): \eta\in \theta_{\CC^m},
\zeta\in F^*(\Der(-\log f))\}}=T^1_{\K_f}F.$$
In the case of $C_\bu(\phi_V)$, we can choose $\ti L_\bu$ to take the form
$$\cd \to \ti L_2\st{d_2}{\too} 
\theta_{\CC^n}\oplus\OO_{\CC^n}\to\OO_{\CC^n}\to 0,$$
with the image of $d_2$ consisting of pairs $(\zeta, b)$ such that
$tf(\zeta)+b\cdot f=0$, so that $\zeta\in \Der(-\log V)$. 
Thus
$$Z_1(C_\bu(\phi_V))=
\{(a,\eta,b)\in \OO\oplus F^*(\theta_{\CC^n})\oplus\OO
: a=t(f\circ F)(\eta)+(f\circ F)b\}.$$
This projects isomorphically to 
$F^*(\theta_{\CC^n})\oplus\OO$ by forgetting the first component. And
$B_1(C_\bu(\phi_V))$ consists of sums of terms of the form
$(0,\zeta,b)\in \OO\oplus F^*(\theta_{\CC^n})\oplus\OO$ such that 
$\zeta\in F^*(\Der(-\log V))$ and $F^*(tf)(\zeta)+b\cdot f\circ F=0$,
coming from $F^*(\ti L_2)$, together with terms of the form
$(--, tF(\eta),a)$ coming from $\ti K_1$ (we do not care what is in the first
component). By projecting this into $F^*(\theta_{\CC^n})\oplus\OO$, 
forgetting the first
component, we see that
$$H_1(C_\bu(\phi_V))=\frac{F^*(\theta_{\CC^n})\oplus\OO}
{\{(\zeta-tF(\eta),b-a):\, F^*(tf)(\zeta)+b\cdot f\circ F=0, \eta\in 
\theta_{\CC^m}, a\in\OO\}}.$$
This is isomorphic to 
$$\frac{F^*(\theta_{\CC^n})}{tF(\theta_{\CC^m})+F^*(\Der(-\log V))},$$
i.e. to $T^1_{\K_V}F$.\eop\\
If $C_\bu(\phi)$ is the cone on a map of complexes $\phi:A_\bu\to B_\bu$, 
there is a long exact sequence of homology
$$\cd \to H_k(A_\bu)\to H_k(B_\bu)\to H_k(C_\bu(\phi))\to H_{k-1}(A_\bu)\to
\cd, $$ constructed by completely standard means (see e.g. \cite{gelman}). 
Thus from \ref{main} we deduce
\begin{corollary}\label{les} There are exact sequences

\begin{equation}\label{seqf}
\cd\to H_1(\k_\bu(f\circ F))
\to
H_1(F^*(\L_\bu))\to T^1_{\K_f}F\to \frac{\OO_{\CC^m}}{J_{f\circ F}}
\to \frac{\OO_{\CC^m}}{F^*(J_f)}\to 0
\end{equation}
and
\begin{equation}\label{seqV}
\cd\to H_1(\ti \k_\bu(f\circ F))
\to
H_1(F^*(\ti \L_\bu))\to T^1_{\K_V}F\to \frac{\OO_{\CC^m}}
{(f\circ F)+J_{f\circ F}}
\to \frac{\OO_{\CC^m}}{F^*((f)+J_f)}\to 0
\end{equation}
in which the maps $T^1_{\K_f}F\to \OO_{\CC^m}/J_{f\circ F}$ and
$T^1_{\K_V}F\to \OO_{\CC^m}/(f\circ F)+J_{f\circ F}$ 
are induced by 
$F^*(tf):\theta(F)\to \OO_{\CC^m}$. \eop
\end{corollary}
\begin{remark}\label{interp}{\em $\OO_{\CC^m}/J_{f\circ F}$ is the $T^1$ of 
$f\circ F$ for right-equivalence. And the morphism 
$T^1_{\K_f}F\to T^1(f\circ F)$ in (\ref{seqf}) 
is a map between deformation functors,
telling us which of the first order deformations of $f\circ F$ we can get
by deforming $F$ alone. A similar statement holds for $T^1_{\K_V}F\to
\OO_{\CC^m}/((f\circ F)+J_{f\circ F})$, with contact-equivalence in place of 
right-equivalence.
}
\end{remark}
In most of what follows we use (\ref{seqf}) to compute $\tau_{\K_f}F:=
\dim_{\CC}T^1_{\K_f}F$ in cases where $f\circ F$ has isolated singularity. In such cases
the Koszul complex $\k_\bu(f\circ F)$ is acyclic. The homology of 
$F^*(\L_\bu)$ computes 
$\mbox{Tor}^{\OO_{\CC^n}}_\bu(\OO_{\CC^n}/J_f,\OO_{\CC^m})$, and 
denoting the rank of the $j$'th Tor module by $\beta_j$, from the exact 
sequence (\ref{seqf}) we obtain
\begin{theorem}\label{betas} If $f\circ F$ has isolated singularity then 
$\tau_{\K_f}F=\mu(f\circ F)-\beta_0+\beta_1$.\eop
\end{theorem}
This equality is the key to the comparisons between $\tau_f(F)$ and the rank
of the vanishing homology of $f\circ F$ under deformations of $F$ alone,
which occupy Section \ref{topdef}.

Note that acyclicity of the Koszul complex $\k_\bu(f\circ F)$ means that for
$k\geq 2, H_k(C_\bu(\phi_f))\simeq 
\Tor_k^{\OO_{\CC^n}}(\OO_{\CC^m},\OO_{\CC^n,0}/J_f)$.

\subsection{Is this the key to any door?}\label{clef}
Where $f\circ F$ has isolated
singularity, by taking the alternating sum of the lengths of the modules
in the exact 
sequences (\ref{seqf}) and (\ref{seqV}) we obtain the formulae 
\be\label{keyf}
\chi(C_\bu(\phi_f))=\mu(f\circ F)-\chi(\OO_{\CC^m,0},\OO_{\CC^n}/J_f)
\ee
and 
\be\label{keyV}
\chi(C_\bu(\phi_V))=\tau(f\circ F)-\chi\bigl(\OO_{\CC^m,0}, 
\OO_{\CC^n}/(f)+J_f\bigr).
\ee
Here, the last term on the right is Serre's intersection
multiplicity (\cite{serre}). This is defined, for modules $M,N$ over the 
ring $R$ ($=\OO_{\CC^n}$ in our case) by
$$\chi(M, N)=\sum_j(-1)^j\ell\bigl(\Tor_j^R(M,N)\bigr).$$
The right hand side makes sense only if $M\otimes_RN$ has finite length,
and in fact this is a sufficient condition for finiteness of all the other
summands.

When $\dim\,M+\dim\,N<\dim\,R$  then $\chi(M,N)=0$\quad (\cite{serre}). 
If also it happens that\linebreak 
$\mbox{Tor}^{\OO_{\CC^n}}_j(\OO_{\CC^m},\OO_{\CC^n}/J_f)=0$ for 
$j>1$ then $\beta_0=\beta_1$ and from \ref{betas} it follows that
$\tau_{\K_f}F=\mu(f\circ F)$. As we shall see in
\ref{surp} below, this explains the surprising equality 
referred to in the opening paragraph of the introduction.

When $f\circ F$ has non-isolated singularity the Koszul complex is no
longer acyclic, and one might wish 
to replace it by a free
resolution $\F_\bu$ of $\OO/J_{f\circ F}$. However, 
in general the comparison morphism 
$\OO/J_{f\circ F}\to\OO/F^*(J_f)$ will not
lift to a morphism of complexes $\F_\bu\to F^*(\L_\bu)$
if $F^*(\L_\bu)$ is not acyclic. In particular,
it is not clear that $tF(\Der(-\log f\circ F))\subset F^*(\Der(-\log f))$,
which is required in order to have such a lift in degree 1. 
The functoriality of the Koszul complex seems to be playing an 
important r\^{o}le here. In order to progress 
towards an understanding of the relation between $\tau_f(F)$ and
the vanishing homology of $V(f\circ F)$ when $f\circ F$ has 
non-isolated singularity, we will need some
understanding of the first Koszul homology of $\OO_{\CC^m}/J(f\circ F)$. 
\subsection{Almost free divisors}\label{afd}
Some of the most interesting developments in the theory of
sections of hypersurface singularities have concerned sections of
free divisors (see e.g. \cite{D589}-\cite{DM}), 
and here the condition of isolated singularity is very far
from being fulfilled. A singular free divisor in $\CC^n$ has singular subspace
of dimension $n-2$, so that among reduced spaces, free divisors have the biggest
possible singular set. An {\it almost free divisor} is a section of a free
divisor $V=f^{-1}(0)$ by a map $F$ which is logarithmically transverse to $V$
outside the origin (this definition is due to Jim Damon); 
so an almost free divisor too is singular in codimension 1.
Thus both the modules $\OO_{\CC^m,0}/F^*(J_f)$ and $\OO_{\CC^m,0}/J_{f\circ F}$
have $(m-2)$ dimensional support. On the other hand in this case 
$\Tor_j^{\OO_{\CC^n}}(\OO_{\CC^m,0}, \OO_{\CC^n}/J_f)=0$ for $j>0$, and so 
the sequence (\ref{seqf}) reduces to the short exact
sequence
\be\label{sesf}
0\to T^1_{\K_f}F\to\frac{\OO_{\CC^m}}{J_{f\circ F}}\to 
\frac{\OO_{\CC^m}}{F^*(J_f)}\to 0
\ee
together with a collection of isomorphisms 
$H_k(C_\bu(\phi_f))\simeq H_{k-1}(\k_\bu(f\circ F))$ 
for $k\geq 2$. It is interesting to note that (\ref{sesf}) allows us to 
give $T^1_{\K_f}F$
the multiplicative structure of the quotient of two ideals:
\be\label{Fman}
T^1_{\K_f}F\simeq \frac{F^*(J_f)}{J_{f\circ F}}.
\ee

The acyclicity of $F^*(\L_\bu)$ for almost free divisors, 
versus the acyclicity of $\k_\bu(f\circ F)$ for sections with 
isolated singularity, shows that these two cases are at opposite 
corners of the field one might wish to survey.

\section{Singularities of matrix families and their determinants}\label{mdets}
In papers \cite{B}, \cite{BT}, \cite{BGZ} and \cite{Ha}, 
parametrised families of $n\times n$ matrices are 
classified up to co-ordinate changes in the parameter
space and parametrised versions of the natural action of 
$\mbox{Sl}_n(\CC)$ and 
$\mbox{Gl}_n(\CC)$. In this section we show that these equivalence relations
are in fact the same as the relations $\K_f$ and $\K_V$ 
when $f$ is the determinant or Pfaffian function on matrix space. 

First we recall the definition of the equivalence relations. 
There are three cases, 
corresponding to symmetric, skew symmetric and arbitrary square matrices, 
each with two flavours, special and general. We remind the reader that
we use $\OO$ to abbreviate $\OO_{\CC^m}$.
\begin{definition}{\em
\ben 
\item
{\it Symmetric matrices:} Symmetric matrix families 
$$S_1,S_2:(\CC^m,0)\to \Sym_n(\CC)$$ 
are 
{\em $\Sl_n$-symmetric equivalent} if there is a matrix family $A\in\Sl_n(\OO)$
and a germ of biholomorphic diffeomorphism $\psi:(\CC^m,0)\to (\CC^m,0)$
such that 
$$S_2=A^t(S_1\circ\psi)A,$$ 
and {\em $\Gl_n$-symmetric equivalent}
if $A$ is allowed to be in $\Gl_n(\OO)$ rather than $\Sl_n(\OO)$.
\item
For {\em skew symmetric matrices}, {\em $\Sl_n$-} and {\em 
$\Gl_n$-skew-equivalence} are defined by the same formulae.
\item
{\it Arbitrary square matrix} families $M_1,M_2:(\CC^m,0)\to \Mat_n(\CC)$
are {\em $\Sl_n$- equivalent} if there are matrix families $A,B\in \Sl_n(\OO)$
and a germ of biholomorphic diffeomorphism $\psi:(\CC^m,0)\to (\CC^m,0)$
such that $$M_2=A(M_1\circ\psi)B,$$
and {\em $\Gl_n$-equivalent} if $A$ and $B$ are allowed to be in $\Gl_n(\OO)$.
\en
}
\end{definition}
It is an immediate consequence of the definitions that the 
extended tangent-spaces to the special and general orbits are\\
$$T_{ss}S=tS(\theta_{\CC^m})+\{A^tS+SA:A\in\sl_n(\OO)\}$$
$$T_{gs}S=tS(\theta_{\CC^m})+\{A^tS+SA:A\in\gl_n(\OO)\}$$
for symmetric and skew symmetric matrices, and
$$T_{sg}M=tM(\theta_{\CC^m})+\{AM+MB:A,B\in\sl_n(\OO)\}$$
$$T_{gg}M=tM(\theta_{\CC^m})+\{AM+MB:A,B\in\gl_n(\OO)\}$$
for arbitrary square matrices. Here the first letter in the subscript
refers to the flavour, special or general, and the second to the type
of the matrix, symmetric, skew or general. We do not need to distinguish,
in our notation, between the symmetric and skew-symmetric cases, since the
equivalence relation is the same. 

We denote the codimension of these tangent spaces (in $S^*(\sym_n(\CC)),
S^*(\sk_n(\CC))$ and $M^*(\mat_n(\CC))$ respectively) by 
$\tau_{ss},\tau_{gs},\tau_{sg}$ and $\tau_{gg}$.

Let $\det$ denote the determinant function on $\Mat_n(\CC)$ in the general
case, and on $\Sym_n(\CC)$ in the symmetric case, and let $V= \{\det=0\}$
(in which space will be clear from the context). Similarly, let 
$\Pf:\Sk_{n}(\CC)\to\CC$ be the Pfaffian function, and in this context let
$V$ denote its zero-locus. Since $\Pf\equiv 0$ if $n$ is odd, from now
on when we are discussing 
skew-symmetric $n\times n$ matrices, $n$ will be assumed even.

\begin{theorem}\label{eqeq} (i) For a symmetric matrix family $S$,
$$T_{ss}S=T\K_{\det}S,\quad T_{gs}S=T\K_VS.$$
(ii) For a general matrix family $M$,
$$T_{sg}M=T\K_{\det}M, \quad T_{gg}M=T\K_VM.$$
(iii) For a skew-symmetric matrix family $S$,
$$T_{ss}S=T\K_{\sPf}S,\quad T_{gs}S=T\K_VS.$$
\end{theorem}
\Proof In each of the three cases, the first equality, concerning
$\Sl_n$-equivalence, can be read off from the well-known free
resolutions mentioned in the introduction. These are described in
detail after this proof. For now, we show
only that the right hand of each of the three pairs of equalities follows from
the left.

In fact we show it only for symmetric matrices --- the other two cases
are essentially identical. 

By comparing the formulae for the tangent 
spaces, we see that
\be\label{sgs}
T_{gs}S=T_{ss}S+\{\lambda S:\lambda\in\OO\}
\ee
since $\gl_n=\sl_n+\{\lambda\cdot \mbox{id}_n:\lambda\in\CC\}$, 
where $\mbox{id}_n$ is the $n\times n$ identity 
matrix, from which 
$\gl_n(\OO)=\sl_n(\OO)+\{\lambda\cdot \mbox{id}_n:\lambda\in\OO\}$ follows.
On the other hand, because $\det$ is a homogeneous function, there is a 
splitting
$$\Der(-\log V)=\Der(-\log\det)\oplus
\OO_{{{\mbox{\footnotesize Mat}}}_n(\CC)}\cdot\chi_e,$$
where $\chi_e$ is the Euler vector field $\sum_{i,j}x_{ij}\p/\p x_{ij}$.
Hence 
\be\label{gsg}
T\K_VS=T\K_{\det} S+S^*(\OO_{\sMat_n(\CC)}\cdot\chi_e)=
T\K_{\det} S+\{\lambda\cdot S:\lambda\in\OO\}.
\ee
From (\ref{sgs}) and (\ref{gsg}) the equality $T_{gs}S=T\K_V S$ follows.
\eop\\

The equalities of the tangent spaces, together with uniqueness of
solutions of ordinary differential equations, imply that the equivalences
themselves coincide, as stated at the beginning of the section. However
we will only need equality of the tangent spaces in what follows.
\section{Free resolutions}\label{alg}
In this section we describe complexes which give free resolutions 
of the jacobian algebras of $\det:\Mat_n(\CC)\to\CC$, 
$\det:\Sym_n(\CC)\to\CC$, and the Pfaffian function $\Pf:\Sk_n(\CC)\to\CC$
for $n$ even. Surprisingly (to us), the original papers where they
appeared make no mention of partial derivatives and vector fields; 
in each of the
three cases, the jacobian ideal is identified as a purely algebraic object,
namely the ideal $I_{n-1}$ generated by the submaximal minors, in the case of 
$\det$, and the ideal $\Pf_{n-2}$ generated by the $(n-2)\times (n-2)$ 
sub-Pfaffians in the case of $\Pf$.
\subsection{The Gulliksen-Neg\aa rd resolution}
For a family $M \in \Mat_n(\OO)$,
Gulliksen and Neg\aa rd constructed in \cite{gn} a complex which 
is a free resolution of 
$\OO/I_{n-1}(M)$ in case the codimension of the variety $V(I_{n-1}(M))$
in $\CC^m$ is 4 (its greatest possible value). Their complex is
\be\label{gnu}
0\to \OO\st{d_4}{\too} {mat}_n(\OO) \st{d_3}{\too}
{sl}_n(\OO)\oplus {sl}_n(\OO)
\st{d_2}{\too} {mat}_n(\OO)\st{d_1}{\too}\OO\too\OO/I_{n-1}(M)\to 0
\ee
where
\begin{itemize}
\item
$d_1(U)={\rm trace}(M^*U)$ where $M^*$ is the adjugate 
of $M$, i.e. the matrix of signed cofactors, 
\item
$d_2(X,Y)=MX-YM$,
\item
$d_3(Z)=(ZM-\frac{tr(ZM)}{n}\ I_n,\ MZ-\frac{tr(MZ)}{n}\ I_n)$,
\item
$d_4(a)=aM^*$
\item we use lower case $\mat(\OO)$ rather than upper case
$\Mat(\OO)$, because we are thinking of the (pull-back of)
the tangent sheaf on the vector space $\Mat(\CC)$. It is 
of course indistinguishable from $\Mat(\OO)$ as $\OO$-module.
\end{itemize}
Several aspects of this resolution come
unexpectedly to our aid. The first is that the
$i,j$'th signed cofactor $M^*_{ij}$ of a
matrix $M$ is equal to 
the partial derivative of $\det$ with respect to the $i,j$'th entry $M_{ij}$,
so in the generic case (where the entries are the variables),
(\ref{gnu}) is a resolution of $\OO_{\sMat_n(\CC)}/J_{\det}$. Moreover
$${\rm trace\,}(M^*X)=\sum_{ij}M^*_{ij}X_{ij}=\sum_{ij}X_{ij}\frac{\p \det}
{\p M_{ij}}=t(\det)(\sum_{ij}X_{ij}\frac{\p}{\p M_{ij}}).$$

The second is that thanks to this fact, we can
interpret the module of relations among the sub-maximal minors as the module
of vector fields annihilating the function $\det$. And the 
Gulliksen-Neg\aa rd  resolution shows that the relations
among them are {\it 
precisely those induced by the action of $\sl_n$}. This fact, which was already
noted by Bill Bruce in \cite{B}, is precisely what is needed to prove
Theorem \ref{eqeq} for matrix families.

Moreover, the Gulliksen-Neg\aa rd complex $GN_\bu(S)$
is equal to $M^*(GN_\bu)$, where $GN_\bu$ is the generic complex, over
$\OO_{\sMat(\CC)}$. Thus it can play the r\^{o}le of the complex $F^*(L_\bu)$
of Theorem \ref{main}.  A similar remark holds for the
other two complexes we now describe.

Once again we have slightly modified the description of the complex
from that found, for example, in \cite{bv}, in order to adapt it to
our situation. 
\subsection{J\'ozefiak's resolution}\label{Sjr}
In \cite{J},
J\'ozefiak constructed a complex of free $\OO$-modules which gives a
resolution of $\OO/I_{n-1}(S)$ provided that the codimension of
the zero variety $V(I_{n-1}(S))$ of $I_{n-1}(S)$ in $\CC^m$
is 3 (its greatest possible value). His complex is
\be\label{j}
0\to \sk_n(\OO)\stackrel{d_3}{\too}
\sl_n(\OO)\stackrel{d_2}{\too}
\sym_n(\OO)\stackrel{d_1}{\too}\OO\to\OO/I_{n-1}(S)\to 0.
\ee
Here
\begin{itemize}
\item
$\sk_n(\OO)$ is the space of order $n$ skew-symmetric matrices
over $\OO$,
\item
$d_1(X)={\rm trace\,}(S^*X)$, where $S^*$ is the adjugate matrix
of $S$,
\item
$d_2(Y)=SY+Y^TS$, and
\item
$d_3(Z)=ZS$.
\end{itemize}
In fact in \cite{J}, $\gl_n(\OO)/\mbox{sk}_n(\OO)$ appears in place of
$\sym_n(\OO)$, and $d_2$ has a different (equivalent) description.

As with the Gulliksen-Neg\aa rd resolution, in the generic case
$d_1$ is equal to $t(\det):\theta_{\sSym(\CC)}\to\OO_{\sSym(\CC)}$, 
and so the acyclicity of (\ref{j}) shows that the module of 
vector fields annihilating $\det$ is generated by the infinitesimal
$\sl_n$ action.
\subsection{The J\'ozefiak-Pragacz resolution}
A complex giving a free resolution of $\OO/\Pf_{n-2}(S)$ 
in case $V(\Pf_{n-2}(S))$  has codimension 6 (its greatest possible
value) in $\CC^m$ is due
to J\'ozefiak and Pragacz \cite{JP}. In slightly modified form it is
$$
0 \to \OO \st{d_6}\too \sk_n(\OO) \st{d_5}\too \sl_n(\OO) \st{d_4}\too
\sym_n(\OO) \oplus \sym_n(\OO) \st{d_3}\too$$
\be\label{jp}
\sl_n(\OO) \st{d_2}\too
\sk_n(\OO) \st{d_1}\too \OO \to \OO/\Pf_{n-2}(S) \to 0
\ee
where
\begin{itemize}
\item
$d_1(U) =\frac{1}{2} \mbox{trace}(S^*U)$, 
where $S^*$ is the matrix of order $(n-2)$
signed Pfaffians of $S$, which satisfies $S^*S= SS^*=\Pf(S)I_n$, 
\item
$d_2(V) = SV + V^TS$,
\item
$d_3(W,X) = S^*W - XS$,
\item
$d_4(Y) = \Big(SY + (SY)^T,\  YS^* + (YS^*)^T \Big)$
\item
$d_5(Z) = ZS - \frac{tr(ZS)}n \ I_n$,
\item
$d_6(a) = aS^*$.
\end{itemize}
In the generic case $\p\Pf/\p S_{ij}=S^*_{ij}$, so 
$\Pf_{n-2}=J_{\sPf}$, 
$d_1=t(\Pf)$, 
(\ref{jp}) is a resolution of the jacobian 
algebra of $\Pf$, and
by acyclicity of (\ref{jp}) it follows that
the vector fields annihilating $\Pf$ are generated by the $\sl_{n}$ action.
\subsection{The morphisms $\phi_f$ for matrix families}
In the case of matrix families, the 
morphism of complexes
$\k_\bu(f\circ F)\st{\phi_f}{\too}F^*(\L_\bu)$ constructed in Section \ref{def}
embodies some non-trivial rules of differentiation. 
For a symmetric or skew-symmetric $n\times n$ matrix family $S$,
$$(\phi_f)_2:\bw^2\theta_{\CC^m}\to \sl_n(\OO)$$
is given by
\be\label{lisy2}
 (\p/\p {x_i} \wedge \p/\p {x_j})\mapsto 
(S^*_{x_i}S_{x_j}-S^*_{x_j}S_{x_i})/2\,
\ee
(where the subscript $x_i$ indicates partial derivative),
and for a  general $n\times n$ matrix family $M$,
$$(\phi_f)_2:\bw^2\theta_{\CC^m}\to \sl_n(\OO)\oplus\sl_n(\OO)$$
is given by
\be\label{lim2}
 (\p/\p {x_i} \wedge \p/\p {x_j})\mapsto 
(M_{x_j}M^*_{x_i}-M_{x_i}M^*_{x_j},  M^*_{x_j}M_{x_i}-M^*_{x_i}M_{x_j})/2\,
\ee

It would be interesting to obtain explicit  formulae for the remaining
$(\phi_f)_j$.  
\subsection{Cohen-Macaulay and Gorenstein properties of determinantal varieties}
\label{cmg}
In each of the three cases, let $\m$ denote the length of the 
resolution (\ref{j}), (\ref{gnu}) and
(\ref{jp}). Since $\m$ is also the codimension of $V(I_{n-1}(S))$,
$V(I_{n-1}(M))$ and $V(\Pf_{n-2}(S))$ in the generic case, 
the three resolutions show that these varieties are Cohen-Macaulay, and the
same conclusion also holds for any matrix family provided the codimension
of the respective variety is $\m$, which is its maximal possible value. 

If $\SS$ is a deformation over base $B$
of a 3-parameter matrix family $S$ meeting this requirement in 
(say) the symmetric case, then the codimension in $\CC^3\times B$ 
of $V(I_{n-1}(\SS))$ is also 3. Since $V(I_{n-1}(S))$ is the fibre of
$V(I_{n-1}(\SS))$ over $0\in B$, $V(I_{n-1}(\SS))$ is finite over $B$, and 
therefore $\OO_{\CC^3\times B}/I_{n-1}(\SS)$ is $\OO_B$-free. This implies
that $\dim_{\CC}\,\OO_{\CC^3,0}/I_{n-1}(S)$ is conserved in a deformation:
it is equal to $\sum_{x}\dim_{\CC}\,\OO_{\CC^3,x}/I_{n-1}(S_t)$ (where the sum
is over the points $x$ into which the isolated zero of $I_{n-1}(S)$ splits), 
since both are equal to the rank of the free sheaf 
$\pi_*\,\OO_{\CC^3\times B}/I_{n-1}(\SS)$. It is the index of intersection of 
the image of $\CC^3$ under $S$
with the set of symmetric matrices of corank greater than 1.  

Similar arguments prove conservation of multiplicity in the other two cases,
when $m=\m$.

Both (\ref{gnu}) and (\ref{jp}) are self-dual complexes, and so if the 
codimension of 
$V(I_{n-1}(M))$ or $V(\Pf_{n-2}(S))$ is $\m$, then 
both are Gorenstein varieties. This has a consequence for the relation between
the Betti numbers of $\GN_\bu(M)$ and $\JP_\bu(S)$ when $m$ is less than 
$\m$, which we now explain.
\begin{lemma}\label{gorlem1} 
Let $R$ be a noetherian local ring, let $\F_\bu=:0\to F_M\to \cd \to F_0\to 0$
be a finite complex of free $R$-modules, and let $\F^\bu$ be its $R$-dual. 
There is a spectral sequence with
$E_2^{p,q}=\Ext^q_R(H_p(\F_\bu),R)$ converging to $H^{p+q}(\F^\bu)$.
\end{lemma}
\Proof Form the double complex $\mbox{Hom}(\F_\bu, I^\bu)$, where $I^\bu$ is
an injective resolution of $R$. Writing the arrow of $\F_\bu$ 
pointing from left to right and the arrow of the injective 
resolution pointing up, the double complex has arrows pointing up
and from right to left.

Taking first the horizontal differential, we get 
$$H_p(\mbox{Hom}(F_\bu, I^q))$$
which is equal to 
$$\mbox{Hom}(H_p(F_\bu),I^q)$$
by injectivity of $I^q$. Now taking vertical differential we get
$$E^{p,q}_2=\mbox{Ext}^q(H_p(F_\bu),R).$$
If we take first the vertical differential in the double complex we get
$$\mbox{Ext}^q(F_p,R)$$
which is equal to $F^p$
for $q=0$ and is zero otherwise (since $\mbox{Ext}$ can also be calculated
by using a projective resolution of $F_p$). Now taking horizontal differential,
we get
$$H^p(\F^\bu).$$
Since both means of calculating the homology of the double complex
must give the same answer, we conclude that the first spectral sequence
must also converge to $H^p(\F^\bu)$.\eop
\begin{lemma}\label{gorlem2} 
Suppose in addition that each of the homology modules
$H_p(\F_\bu)$ has finite length and that $R$ is a regular local
ring of dimension $m$. Then
$$H^p(\F^\bu)\simeq \Ext^m(H_{p-m}(\F_\bu),R).$$
\end{lemma}
\Proof Each module $H_p(\F_\bu)$ now has a free resolution of length $m$.
By the Lemma of
Ischebeck (p. 133 of \cite{matsu}) 
$$\mbox{Ext}^q(H_p(F_\bu),\OO)=0$$
except when $q=m$. It follows that the spectral sequence of
\ref{gorlem1} collapses at $E_2$, and the lemma follows.\eop\\
Now suppose that the complex $\F_\bu$ is self-dual, in the sense that
there is an integer $\m$ (the length of the complex) 
such that $H^p(F^\bu)=H_{\m-p}(F_\bu)$. This is the
case for the complexes $\GN_\bu(M)$ and $\JP_\bu(S)$, since they 
are the pull-backs of
free resolutions of Gorenstein quotients. Then \ref{gorlem2} gives
\be\label{goreq1}
H_{\m-p}(\F_\bu)\simeq \Ext^m(H_{p-m}(\F_\bu),R).
\ee
\begin{proposition}\label{gorpl}
Suppose that $\F_\bu$ is a self-dual complex of 
free modules, of length $\m$, over the regular local 
ring $R$ of dimension $m$,
with all homology modules of $\F_\bu$ having finite length.  
Then 
$H_{\m-p}(F_\bu)$ has the same length as $H_{p-m}(\F_\bu).$
\end{proposition}
\Proof Immediate from (\ref{goreq1}) and the following lemma.\eop
\begin{lemma}\label{gorlem3} 
Let $M$ be an $R$-module of finite length. Then 
$\Ext^m(M,R)$ has the same length as $M$.
\end{lemma}
\Proof 
Because $R$ is Gorenstein, local duality (see e.g. \cite{BH} Chapter) 
gives us
$$\Ext^m(M,R)=\Ext^m(M,\omega_R)
=\mbox{Hom}(H^0_{\{0\}}(M),E(k)),$$
where $\omega_R$ is the dualising module of $R$ and 
$E(k)$ is the injective hull of the residue field of $R$.
Since $M$ is supported only at $0$, this gives
$$\mbox{Ext}^m(M,R)=\mbox{Hom}(M,E(k)).$$
Finally, an easy induction (\cite{BH}, page 102) shows that the
length of $\mbox{Hom}(M,E(k))$ is equal to the length of $M$ for any
$R$-module $M$ of finite length.
\eop 
\begin{proposition}\label{gorp} (i) Let $M$ be a general matrix family on
$m<4$ parameters, and suppose that $\det(M)$ has isolated singularity.
Then $$\beta_k(\GN_\bu(M))=\beta_{k+(4-m)}(\GN_\bu(M)).$$
(ii) Let $S$ be a skew-symmetric matrix family on $m<6$ parameters, and
suppose that $\Pf(S)$ has isolated singularity. Then
$$\beta_k(\JP_\bu(S))=\beta_{k+(6-m)}(\JP_\bu(S)).$$ 
\end{proposition}
\section{$\tau, \mu$ and the vanishing homology of sections with isolated 
singularity}\label{topdef}
We will suppose, throughout this section, that $\OO_{\CC^n}/J_f$ is 
Cohen-Macaulay, of dimension $n-\m$. Then
the dimension of $\OO_{\CC^m}/F^*(J_f)$ is at least $m-\m$, and if
it is $m-\m$ then $\OO_{\CC^m}/F^*(J_f)$ is Cohen-Macaulay. 
Moreover in this case if $\L_\bu$ is a free $\OO_{\CC^n}$-resolution 
of $\OO_{\CC^n}/J_f$ then
$F^*(\L_\bu)$ is a free $\OO_{\CC^m}$-resolution of $\OO_{\CC^m}/F^*(J_f)$. 
From Theorem \ref{betas} we therefore 
obtain
\begin{theorem}\label{imax} If $m=\m$ and $f\circ F$ has isolated singularity
then 
$$\tau_{\K_f}F=\mu(f\circ F)-\dim_{\CC}\OO_{\CC^m}/F^*(J_f).$$\eop
\end{theorem} 
Applying this to the matrix families we considered in Section \ref{mdets}, we
have
\begin{corollary}\label{max} 
(i) Symmetric case with $m=3$: 
$$\tau_{ss}(S)=\mu(\det(S))-\dim_{\CC}\OO_{\CC^3}/I_{n-1}(S).$$
(ii) General case with $m=4$:
$$\tau_{sg}(M)=\mu(\det(M))-\dim_{\CC}\OO_{\CC^4}/I_{n-1}(M).$$
(iii) Skew-symmetric case with $m=6$:
$$\tau_{ss}(S)=\mu(\Pf(S))-\dim_{\CC}\OO_{\CC^6}/\Pf_{n-2}(S).$$
\end{corollary}
If $m=\m -1$ or $m=\m -2$ a surprising phenomenon occurs:
\begin{lemma}\label{surp} (i) 
Suppose $m=\m -1$, and $f\circ F$ has isolated 
singularity. Then 
the numbers $\b_0$ and $\b_1$ in Theorem \ref{betas}
are finite and equal to one another, and $\beta_j=0$ for $j\geq 2$. 
In consequence, $\tau_{\K_f}F=\mu(f\circ F)$.\\
(ii) If $m=\m-2$ and $f\circ F$ has isolated singularity, then
$\b_0$, $\b_1$ and $\b_2$ are finite and
$\b_j=0$ for $j\geq 3$. Moreover 
$\b_0+\b_2=\b_1$, so that 
$\tau_{\K_f}F=\mu(f\circ F)+\beta_2$.\\
\end{lemma}
\Proof We have already remarked in Subsection \ref{clef} that given
the vanishing of higher Tor modules, the relations between the $\beta_j$
follow from the vanishing of $\chi(\OO_{\CC^n}/J_f,\OO_{\CC^m})$. Since we
need this vanishing, however, we give a self-contained proof.

(i) Choose a deformation $\s F$ of $F$ on one extra parameter $t$, so 
that $\OO_{\CC^{\m-1}\times\CC}/{\s F}^*(J_f)$ has dimension $0$. 
Then $\s F^*(\L_\bu)$ is acyclic. The short exact sequence of complexes
$$0\to \s F^*(\L_\bu)\st{t\cdot}{\too}\s F^*(\L_\bu)\to F^*(\L_\bu)\to 0$$
gives rise to a long exact sequence of homology; since $\s F^*(\L_\bu)$
is acyclic, this ends
$$0=H_1(\s F^*(\L_\bu))\to H_1(F^*(\L_\bu))\to 
H_0(\s F^*(\L_\bu))\to
H_0(\s F^*(\L_\bu))\to
H_0(F^*(\L_\bu))\to 0.$$
As $H_0(\s F^*(\L_\bu))=\OO_{\CC^{\m-1}\times\CC}/\s F^*(J_f)$,
it has finite length, and the conclusion follows from the fact that
the alternating sum of the lengths of the modules in an exact sequence is 0.\\
(ii) Let $\s F_1$ be a deformation of $F$ on the parameter $t_1$
and let $\s F_2$ be a deformation of $\s F_1$ on the parameter $t_2$,
such that $V(\s F_2^*(J_f))$ has codimension $\m$. The argument of (i)
applied to the long exact sequence arising from
$$0\to \s F_2^*(\L_\bu)\st{t_2\cdot}{\too}\s F_2^*(\L_\bu)\to\s F_1^*(\L_\bu)\to 0$$
shows that $\beta_2(\s F_1^*(\L_\bu))=0$ and 
$\beta_1(\s F_1^*(\L_\bu))<\infty$. An analogous argument, applied to the 
(longer) exact sequence arising from 
$$0\to \s F_1^*(\L_\bu)\st{t_1\cdot}{\too}\s F_1^*(\L_\bu)\to F^*(\L_\bu)\to 0$$
then gives the result.
\eop\\
It is curious in case (i) that despite the equality of $\mu$ and $\tau$, 
the natural map from 
$T^1_{\K_f}F$ to the jacobian algebra is not an isomorphism.

For our matrix families we conclude from \ref{betas} and \ref{surp}(i)
\begin{corollary}\label{submax} 
(i) Symmetric case with $m=2$: 
$$\tau_{ss}(S)=\mu(\det(S)).$$
(ii) General case with $m=3$:
$$\tau_{sg}(M)=\mu(\det(M)).$$
(iii) Skew-symmetric case with $m=5$:
$$\tau_{ss}(S)=\mu(\Pf(S)).$$
\eop\end{corollary}
Also somewhat surprisingly, under a mild assumption on the distributions
$\Der(-\log V)$ and $\Der(-\log f)$, 
the two disparate phenomena described in
\ref{imax} and \ref{surp} are both subsumed into the same phenomenon
when we consider the vanishing homology of $F^{-1}(V)$ under deformation
of $F$.

The assumption is the following:
\ben
\item
There exist perturbations of $F$ which are logarithmically transverse
to $V$, and
\item
at each point $x\in V$ where $\dim_{\CC}T^{\log}_xV\geq n-m,$ 
\be\label{eq}
 T^{\log}_xV=\Der(-\log f)(x). 
\ee
\en
This equality holds, for example, if $f\in m_xJ_f$. 
For if $f\in  J_f$ then $\Der(-\log V)$ splits as a direct
sum $\Der(-\log f)\oplus\OO_{\CC^n}\cdot\chi$ where $\chi$ is a vector
field such that $\chi\cdot f = f$. So if in addition we can choose $\chi$ to
vanish at $x$ (i.e. if
$f\in m_xJ_f$) then (\ref{eq}) follows. 
In fact if $f\in m_xJ_f$ and $u$ is any unit then 
$uf\in m_xJ_{uf}$, so (\ref{eq})
holds for {\it every} choice of equation.
This has a partial converse:
\begin{proposition} Suppose that $T^{\log}_xV\neq 0$. Then (17) holds 
for every choice of equation for $V$ if and only if $f\in m_xJ_f$.
\end{proposition}
\Proof ``If'' is already shown. 
Choose $\chi\in\Der(-\log f)$ such that $\chi(x)
\neq 0$. If $\chi(x)\in \Der(-\log uf)(x)$, then there
must exist $\eta\in m_x\theta_{\CC^n,x}$ such that $(\chi+\eta)\cdot uf=0$.
This gives
$$u (\eta\cdot f)+f\bigl((\chi+\eta)\cdot u\bigr)=0.$$
Now choose a unit $u$ such that $d_xu(\chi(x))\neq 0$. Then 
$(\chi+\eta)\cdot u$ is a unit in $\OO_{\CC^n,x}$, so $f\in m_xJ_f$.\eop\\

We will say that $m$ is {\it in the range of holonomy} with respect to $f$ if 
(\ref{eq}) holds at all points of $V$ where $\dim_{\CC}T^{\log}_xV\geq n-m$. 
\begin{theorem}\label{smax}  Suppose $m=\m$ or $m=\m-1$, and is in 
the range of holonomy
with respect to $f$, and let $F_t$ be a 
perturbation of $F$ which is logarithmically transverse to $V$. 
Let $X_t=V(f\circ F_t)$
and $X_0=V(f\circ F)$.
Then $X_t$ has the homotopy type of a wedge of 
$\tau_{\K_f}F$ copies of the $(m-1)$-sphere. 
\end{theorem}
\Proof When $m=\m -1$ the argument is straightforward: $X_t$ is a Milnor fibre
of $X_0$, since $F$ can only meet $V$ at regular points, 
and there transversely.
Hence $X_t$ has the homotopy type of a wedge of $\mu$ spheres. By
\ref{surp}(i), $\mu=\tau$.

When $m=\m$, $X_t$ is no longer a smoothing of $X_0$. Instead, it has 
an isolated singular point at each zero of $F_t^*(J_f)$.
However, $T^1_{\K_V}F_t$ 
is everywhere 0, so by the assumption that $m$ is in the range of holonomy
of $f$, $T^1_{\K_f}F_t$ vanishes also at each point of $V(F_t^*(J_f))$ (these 
points all lie on $F_t^{-1}(V)$), and so by \ref{betas}, the Milnor 
number of 
$f\circ F_t$ at each singular
point $x$ is equal to the local contribution $\beta_0(F_t,x)$. 
By smoothing each singularity, and so obtaining a Milnor fibre for the
isolated singularity $(X_0,0)$, we would 
increase the rank of the middle homology by
the local Milnor number of $X_t$ at $x$.
It follows that the rank of $H_{m-1}(X_t)$ is
$\mu-\sum_x\beta_0(F_t,x)$. However, $\sum_x\b_0(F_t,x)=\b_0$, by the 
well-known argument sketched at the start of Subsection \ref{cmg}. 

It is also 
well known that every fibre of a deformation of an isolated hypersurface 
singularity, whether smooth or not, has the homotopy type of a
wedge of spheres. \eop\\
\begin{remark}\label{normalforms}{\em
(i) The existence of a perturbation $F_t$ of $F$ which is logarithmically
transverse to $V$ is not always assured (see (ii) below); 
however, the logarithmic 
stratification of matrix space is finite, and an argument of Jim Damon's 
(\cite{leg}, 2.4) using Sard's theorem shows that for families of matrices
the required perturbations 
do exist. Moreover, standard row-reduction arguments show that the varieties
$\det=0$ and $\Pf=0$ are everywhere locally quasihomogeneous, so that
the range of holonomy of $\det$ and $\Pf$ has no upper bound. 
In fact,
at the singular points of the determinant of a generic 
matrix in 4 parameters or of a generic symmetric matrix in 3 parameters, 
and at a singular point of the Pfaffian of a
generic skew-symmetric matrix in 6 parameters, the families are 
$\Sl_n$-equivalent,
respectively, to 
\be\label{norm}
\left(
\begin{array}{ccc}
x_1&x_2&0\\
x_3&x_4&0\\
0&0&M_{n-2}
\end{array}
\right),
\left(
\begin{array}{ccc}
x_1&x_2&0\\
x_2&x_2&0\\
0&0&S_{n-2}
\end{array}
\right),
\quad \mbox{and}\quad
\left(
\begin{array}{ccccc}
0&x_1&x_2&x_3&0\\
-x_1&0&x_4&x_5&0\\
-x_2&-x_4&0&x_6&0\\
-x_3&-x_5&-x_6&0&0\\
0&0&0&0&A_{n-4}
\end{array}
\right),
\ee
where $M_{n-2}$, $S_{n-2}$ and $A_{n-4}$ are constant matrices with
non-vanishing determinants.
The determinants and Pfaffian of the families in (\ref{norm}) are equal to 
\be
x_1x_4-x_2x_3,\quad x_1x_3-x_2^2,\quad \mbox{and}\quad x_1x_6-x_2x_5+x_3x_4,
\ee
each of which has a non-degenerate critical point at $0$. The ideals
$I_{n-1}$ and $\Pf_{n-2}$ are in each case equal to the maximal ideal.\\

(ii) Consider the hypersurface $V=\{y(x+y)(x-y)(x+zy)=0\}$ in 
$\CC^3$ and the map $F(x,y)=(x,y,0)$. This is the total
space of a family of quadruple lines in the plane, with parameter $z$. 
As the cross-ratio varies with $z$, $V$ is not analytically trivial along the
$z$-axis.
We claim that 
$F$ has no perturbation which
is logarithmically transverse to $V$. For on the one hand 
at every point $P$ on the $z$-axis,
$T^{\log}_PV=0$, so that no map from $\CC^2$ to $\CC^3$ can meet $V$
transversely at $P$, while on the other hand since $F$ meets the $z$-axis
in an isolated point,
every perturbation $F_t$ of $F$ will also meet the $z$-axis. 
In this example, $V$ is neither globally weighted homogeneous, 
nor locally quasihomogeneous
at any point on the $z$-axis.\\

(iii) In contrast, the hypersurface (cf \cite{cmn}) with equation 
$$f(x,y,z)=x^5z+x^3y^3+y^5z,$$
is globally homogeneous, 
but not locally quasihomogeneous at any point of
the $z$-axis outside $0$. As a {\it Macaulay} calculation readily shows,
$\Der(-\log f)$ is generated by vector fields which vanish everywhere on the
$z$-axis, so $\Der(-\log f)(0,0,z)=0$. On the other hand 
$\Der(-\log V)(0,0,z) \neq 0$ for $z\neq 0$, since it contains the value of
the Euler vector field. The section $F(x,y)=(x,y,x+y)$, with 
$\tau_{\K_f(F)}= 10$, has a perturbation 
$F_t(x,y)=(x,y,x+y+t)$ which is logarithmically transverse to $V$ ---
in fact it is transverse to the distribution spanned by the Euler vector field
--- but nevertheless we have 
$F_t(0,0)\in V$ and $(T^1_{\K_f}F_t)_{(0,0)}\neq 0$. 

We must admit that here $\OO_{\CC^3}/J_f$ is not Cohen-Macaulay, so
that the conclusion of Theorem \ref{smax} fails for other reasons too.
Indeed the conclusion of Lemma \ref{imax}(i) fails also: one calculates that 
$\mu(f\circ F)=25$ and $\dim_{\CC}\OO_{\CC^2}/F^*(J_f)=19$, so
that here
$$\tau_{\K_f}(F)\neq\mu(f\circ F)-\dim_{\CC}\OO_{\CC^3}/F^*(J_f).$$ 
}

\end{remark}
From \ref{smax} and \ref{normalforms}(i) we conclude:
\begin{theorem} In each of the three kinds of matrix families, if $m=\m$
or $m=\m -1$ then the relevant Tjurina number $\tau$ is equal to the
rank of the vanishing homology of the determinant or Pfaffian.\eop
\end{theorem}
One 
further comparison between $\tau_{\K_f}(F)$ and the vanishing homology
follows immediately from \ref{betas}, \ref{gorpl} and \ref{surp}(ii) 
in case $\OO_{\CC^n}/J_f$ is Gorenstein.
\begin{proposition}\label{gawp} 
Suppose that $\OO_{\CC^n}/J_f$ is Gorenstein of dimension $n-\m$, and that
\linebreak
$F:(\CC^{\m-2},0)\to (\CC^n,0)$ has $\tau_{\K_f}F<\infty$. Then
$$\tau_{\K_f}F=\mu(f\circ F)+\dim_{\CC}\,\OO_{\CC^{\m-2}}/F^*(J_f).$$
\eop
\end{proposition}
\begin{corollary}\label{ssmax} (i) If $M$ is a general matrix family on two
parameters, then 
$$\tau_{sg}(M)=\mu(\det(M))+\dim_{\CC}\,\OO/I_{n-1}(M).$$
(ii) If $S$ is a skew-symmetric matrix family on four parameters, then
$$\tau_{ss}(S)=\mu(\Pf(S))+\dim_{\CC}\,\OO/\Pf_{n-2}(S).$$
\eop
\end{corollary}
\begin{remark}{\em
Since $\J_\bu$ is not self-dual, an analogue of \ref{ssmax} does not hold
in general for symmetric matrices. According to \cite{B},
any 1-parameter family of such
matrices with finite Tjurina number is equivalent to 
$S=\mbox{diag}\{x^{a_1}, x^{a_2}, \dots, x^{a_n}\}$, for some non-decreasing 
sequence of integers $0 \le a_1 \le a_2 \le \dots \le a_n$. For such a family
$S$,
$$\tau_{ss}(S) = \sum_{i=1}^n (n-i+1)a_i -1\,, \qquad
\mu(det\,S) = \sum_{i=1}^n a_i - 1\,, \qquad 
\b_0(\JJ_\bu) = \sum_{i=1}^{n-1} a_i\,.
$$
Therefore, $\tau_{ss}(S) = \mu(det\,S) + \b_0$ if and only if
the corank of the matrix $S(0)$ is at most 2, in which case $S^*(J_{\det})$
is actually a complete intersection ideal
and part 2 of Theorem \ref{Ck} (below) applies.} 
\end{remark}
\subsection{Sections of isolated hypersurface singularities}
\label{Sk}
Now suppose that $f:(\CC^n,0)\to (\CC,0)$ 
has an isolated singularity at the origin. Then $\m=n$, and
there are three values of $m$ when the number $\b_1$ in
the right-hand side of (\ref{betas}) is easy to calculate:
$\b_1 = 0, \  \b_0, \  2\b_0$ if $m=n,\  n-1, \  n-2$  
respectively. From \ref{imax}, \ref{surp} and \ref{gawp}
follows

\begin{theorem}\label{Ck} {\em \cite{BGZ}\\

0) $\tau_{\K_f}(F) = \mu(f \circ F) - \b_0$\hspace{.7in} if $m=n$;

1) $\tau_{\K_f}(F) = \mu(f \circ F)$ \hspace{1in} if $m=n-1$;

2) $\tau_{\K_f}(F) = \mu(f \circ F) + \b_0$ \hspace{.7in}if $m=n-2$.}\eop
\end{theorem}
If $m=n-1$, a generic perturbation $F_t$ of $F$ will 
be transverse to $V(f)$
and miss $0\in\CC^n$ altogether. 
Hence $V(f\circ F_t)$ is a smoothing of
$V(f\circ F)$.
If $m=n$, $\b_0 = \mu(f) \cdot {\rm deg\,}F$; moreover
if $F_t$ is a generic perturbation of $F$, it will
cover $0\in\CC^n$ $\deg\,F$ times, and thus $F_t^{-1}(V(f))$
will have $\deg\,F$ singular points, each with Milnor number $\mu(f)$.
Smoothing these to get a Milnor fibre of $f\circ F$, we 
increase the rank of the middle homology of
$F_t^{-1}(V(f))$ by $\deg\ F\cdot\mu(f)$; hence 
\begin{corollary}If $m=n$ or $m=n-1$, 
$$\mbox{rank}\ H_{m-1}(V(f\circ F_t))= \tau_{\K_f}(F).$$
\eop
\end{corollary}
\section{Cohen-Macaulay properties of the relative $T^1$}\label{CM}
In the last section we showed that given a diagram 
$$(\CC^m,0)\st{F}{\too}(\CC^n,0)\st{f}{\too}(\CC,0)$$
with $\OO_{\CC^n}/J_f$ Cohen-Macaulay of codimension $\m$, $\tau_{\K_f}F<\infty$,
and $m=\m$ or $m=\m-1$, then
$$\mbox{rank(vanishing homology of}\ \  V(f\circ F))
=\tau_{\K_f}F.$$
In each case, this formula was proved by placing $T^1_{\K_f}F$
in an exact sequence which related it to other modules whose lengths are
conserved in a deformation. Now suppose that
\linebreak $\f:(\CC^m\times B,(0,0))\to(\CC^n,0)$
is a deformation of $F$ over the smooth base $B$. A slight modification of this 
argument gives us information about the relative $T^1$, 
$$T^1_{\K_f/B}\f:=\frac{\theta(\f)}{t\f(\theta_{\CC^m\times B/B})+
\f^*(\Der(-\log f))}.$$
\begin{theorem}\label{cm}
Suppose that in the diagram above, $\tau_{\K_f}F<\infty$. Then\\
(i) if $m=\m$ then $T^1_{\K_f/B}\f$ is Cohen-Macaulay over $\OO_{\CC^m\times B}$,
of dimension equal to $\dim B$. Moreover, it is free of rank $\tau_{\K_f}F$ over
$\OO_B$.\\
(ii) if $m=\m-1$ and in addition we suppose that 
$\codim\,\OO_{\CC^m\times B}/\f^*(J_f)=\m$, then the same
conclusions hold as in (i).
\end{theorem}
\Proof 
(i) From the relative version of Theorem \ref{main}, we obtain the
exact sequence
$$0\to \mbox{Tor}_1^{\OO_{\CC^n}}(\OO_{\CC^n}/J_f,\OO_{\CC^m\times B})
\to T^1_{\K_f/B}\f\to \OO_{\CC^m\times B}/J^{\mbox{rel}}_{f\circ\f}
\to \OO_{\CC^m\times B}/\f^*(J_f)\to 0.$$
As the absolute $\mbox{Tor}$, 
$\mbox{Tor}_1^{\OO_{\CC^m}}(\OO_{\CC^n}/J_f,\OO_{\CC^m})$, vanishes, so
does the parametrised $\mbox{Tor}$, and this sequence reduces to
a short exact sequence. The second and third modules in this short exact
sequence are finite over $\OO_B$, and are both of dimension equal
to $\dim B$. It follows that they are $\OO_B$-free. Hence, by the depth lemma,
so is the first.\\
(ii) The argument is almost identical. The only differences are that
we have explicitly to require
that $\codim\, V(\F^*(J_f))=\m$, in order to 
guarantee the vanishing of the relative $\mbox{Tor}$, and that instead of 
$\OO_{\CC^m\times B}/\f^*(J_f)$ being free over $\OO_B$, it is finite of
dimension $\dim\,B-1$. Once again, the depth lemma guarantees that 
$ T^1_{\K_f/B}\f$ is $\OO_B$-free.\eop\\

From this one can prove Theorem 4.5 by the argument of Section 5
of \cite{DM}.\\

Recall from Section \ref{def} the Damon 
module $T^1_{\K_V}S$, 
in which $\Der(-\log \det)$ is replaced by $\Der(-\log V)$,
and its relative version $T^1_{\K_V/B}\SS$.
Their support is the set of points where $S$ is not
logarithmically transverse to $V$ (respectively, ${\s S}$ is not 
logarithmically
transverse to $V$ relatively to $B$). 
The {\it discriminant} ${\s D}\subset B$ is defined to be the projection
to $B$ of the support of $T^1_{\K_V/B}\SS$.\\

It tuns out that in all simple singularities of symmetric
matrix families in 2, 3 {\it and 4} variables, ${\s D}$ is a free divisor.
For simple families in 2 variables, this 
follows from the explicit description of the discriminant given in Proposition
3.3 of
\cite{GZ}; as implied also by \cite{B}, for 3 variables
it follows from the fact that all of the 
discriminants are the discriminants of simple singularities of
functions on a manifold with boundary. 
The result in dimension 4 (observed empirically by computer
calculation) is surprising; also surprising in these 
(rather few) examples of 4-parameter simple symmetric matrix families 
is that the conclusion of Theorem \ref{cm} (and therefore Theorem \ref{smax})
holds for them. We note that by the argument of Jim Damon (\cite{leg}),
the freeness of $\s D$ is closely related to the (experimentally verified)
fact that $T^1_{\K_V/B}\s S$ is Cohen-Macaulay of dimension $\dim\,B-1$,
and that Damon's condition on the existence of ``Morse-type singularitiees''
holds.

Freeness of the discriminant in the base of a family of general matrices in
two variables would follow from Conjecture 3.5 of \cite{GZ}.

\end{document}